\documentclass[a4paper,11pt,twoside]{article}

\usepackage{amssymb}
\usepackage[T1]{fontenc}
\usepackage[utf8]{inputenc}
\usepackage{mathrsfs}
\usepackage{amsfonts,amsmath,lmodern,fancyhdr,lastpage,graphicx,nccfoots,caption}
\usepackage{afterpage}
\usepackage{geometry}
\usepackage{amsthm}
\usepackage{amstext}
\usepackage{amssymb, url}

\usepackage{graphicx}

\makeatletter

\usepackage{hyperref}

\vfuzz \hfuzz
\newcounter{a}
\ifodd\thea\else\stepcounter{a}\fi

\fancypagestyle{plain}{%
\fancyhf{}
\fancyfoot[C]{\tiny{Encontro Nacional da SPM 2014, Geometria e Topologia, pp.\ \thea-\pageref{LastPage}}}

}
\numberwithin{equation}{section}
\numberwithin{figure}{section}
\theoremstyle{plain}
\newtheorem{thm}{Theorem}[section]
  \theoremstyle{definition}
  \newtheorem{defn}[thm]{Definition}
  \theoremstyle{remark}
  
  \theoremstyle{plain}
  
  \theoremstyle{plain}
  \newtheorem{prop}[thm]{Proposition}
  \theoremstyle{plain}
  \newtheorem{cor}[thm]{Corollary}

\def\quot{/\!\!/}
\def\hom{\mathsf{Hom}}
\renewcommand{\leq}{\leqslant}

\renewcommand{\geq}{\geqslant}

\newcommand{\bG}{\mathbf{G}}
\newcommand{\X}{\mathfrak{X}}

\newcommand{\R}{\mathbb{R}}

\newcommand{\C}{\mathbb{C}}

\newcommand{\VV}{\mathbb{V}}

\newcommand{\GL}{\mathrm{GL}}

\newcommand{\SL}{\mathrm{SL}}

\newcommand{\la}{\langle}
\newcommand{\ran}{\rangle}

\DeclareMathOperator{\Aut}{Aut}

\DeclareMathOperator{\tr}{tr}

\DeclareMathOperator{\End}{End}

\DeclareMathOperator{\Fix}{Fix}

\newcommand{\liep}{\mathfrak{p}}

\newcommand{\liek}{\mathfrak{k}}

\newcommand{\lieg}{\mathfrak{g}}
\newcommand{\liegc}{\mathfrak{g}^{\mathbb{C}}}

\begin{document}
%\label{pp}
\thispagestyle{plain}

%%%%%%%%%%%%%%%%%%%%
% T’tulo do artigo
\begin{center}
\Large
\textsc{On the homotopy type of free group character varieties\footnote{A.C., C.F. and A.O. were partially supported by FCT projects PTDC/MAT/120411/2010, MAT-GEO/0675/2012 and the Reseacrh Unit PEst-OE/MAT/UI4080/2011, Portugal; S.L. was partially supported by  U.S. NSF grants DMS 1107452, 1107263, 1107367 - Gear Network, and 1309376, and the Simons Foundation grant 245642.}}
\end{center}

%%%%%%%%%%%%%%%%%%%%
% Primeiros autores, com endereo comum

\begin{center}
\textit{Ana Casimiro} \smallskip \\
\begin{tabular}{l}
\small Departamento de Matem\'{a}tica, Faculdade de Ci\^{e}ncias e Tecnologia,\\
\small Universidade Nova de Lisboa, \\
\small Quinta da Torre,  2829-516 Caparica, Portugal       \\
\small e-mail: \texttt{amc@fct.unl.pt}   \\
\end{tabular}
\end{center}

%%%%%%%%%%%%%%%%%%%%
% Outro autor, com endereo diferente

\begin{center}
\textit{Carlos Florentino} \smallskip \\
\begin{tabular}{l}
\small Departamento de Matem\'{a}tica, Instituto Superior T\'{e}cnico,     \\
\small Universidade de Lisboa,                     \\
\small Avenida Rovisco Pais, 1049-001 Lisboa, Portugal       \\
\small e-mail: \texttt{cfloren@math.tecnico.ulisboa.pt}
\end{tabular}
\end{center}

%%%%%%%%%%%%%%%%%%%%
% Outro autor, com endereo diferente

\begin{center}
\textit{Sean Lawton} \smallskip \\
\begin{tabular}{l}
\small Department of Mathematical Sciences, George Mason University                  \\
\small 4400 University Drive, Fairfax, Virginia, 22030, USA       \\
\small e-mail: \texttt{slawton3@gmu.edu}
\end{tabular}
\end{center}

%%%%%%%%%%%%%%%%%%%%
% Outro autor, com endereo diferente

\begin{center}
\textit{Andr\'{e} Oliveira} \smallskip \\
\begin{tabular}{l}
\small Departamento de Matem\'{a}tica,             \\
\small  Universidade de Tr\'{a}s-os-Montes e Alto Douro,                \\
\small Quinta dos Prados, 5001-801 Vila Real, Portugal       \\
\small e-mail: \texttt{agoliv@utad.pt}
\end{tabular}
\end{center}

%%%%%%%%%%%%%%%%%%%%
%% Resumo (em portugus)
%
%\noindent
%\textbf{Resumo:}
%Seja $G$ um grupo alg\'{e}brico, redutivo e real, com subgrupo compacto maximal $K$, e seja $F_{r}$ um grupo livre de $r$ geradores.
%Mostramos que o espa\c{c}o de \'{o}rbitas fechadas em $\hom(F_{r},G)/G$
%admite uma deforma\c{c}\~{a}o retractil forte para o espa\c{c}o de \'{o}rbitas  $\hom(F_{r},K)/K$.
%Em particular, todos estes espa\c{c}os t\^{e}m o mesmo tipo de homotopia. Calculamos os polin\'{o}mios de Poincar\'{e} destes espa\c{c}os, para grupos de dimens\~{a}o pequena, $G$, tais como  $\Sp(4,\mathbb{R})$ e $\U(2,2)$.
%
%
%\medskip

%%%%%%%%%%%%%%%%%%%%
% Abstract (em ingl\^es)

\noindent
\textbf{Abstract}
Let $G$ be a real reductive algebraic group with maximal compact
subgroup $K$, and let $F_{r}$ be a rank $r$ free group.
Here, we summarize the construction of a natural strong deformation retraction
from the space of closed orbits in $\hom(F_{r},G)/G$
 to the orbit space $\hom(F_{r},K)/K$.
In particular, these spaces have the same homotopy type.
 %We compute
%the Poincar\'{e} polynomials of these spaces for some low rank groups
%$G$, such as $\Sp(4,\mathbb{R})$ and $\U(2,2)$.
\medskip

%%%%%%%%%%%%%%%%%%%%%
%% Lista de palavras-chave
%
%\noindent\textbf{palavras-chave:} Variedades de caracteres; Grupos reductivos reais; Variedades de representa\c{c}\~oes.
%
%\medskip

%%%%%%%%%%%%%%%%%%%%
% keywords (vers\~ao inglesa da lista de palavras-chave)

\noindent\textbf{keywords:}
Character varieties; Real reductive groups; Free group representations.

%%%%%%%%%%%%%%%%%%%%
%%%%%%%%%%%%%%%%%%%%
% Primeira sec‹o

\section{Introduction}
In this article, we present one of the main results in
\cite{Casimiro-Florentino-Lawton-Oliveira:2014}, about the
homotopy equivalence between two related moduli spaces of representations
of a free group $F_r$ on $r$ generators.
These are the $G$-character variety $\mathfrak{X}_{r}(G)$, consisting of the closed
orbits in the conjugation quotient space  $\hom(F_{r},G)/G$, where $G$ is a real reductive
group, and the related quotient space $\mathfrak{X}_{r}(K)=\hom(F_{r},K)/K$,
where $K$ is a maximal compact
subgroup of $G$.
One application of this result was the computation of the Poincar\'{e}
polynomials of some of these character varieties, for non-compact $G$.  We also studied the
relation between the topology and \emph{geometry} of
the character varieties $\mathfrak{X}_{r}(G)$ and (the real points
of) $\mathfrak{X}_{r}(\mathbf{G})$, where $\mathbf{G}$ is the complexification of $G$, making explicit use of \emph{trace
coordinates}. We provided a detailed analysis of
some examples (real forms $G$ of $\mathbf{G}=\SL(2,\mathbb{C})$),
showing how the geometry of these spaces compare, and how to understand
the deformation retraction in these coordinates.
We also briefly described the Kempf-Ness sets for some of these examples.
We thank the referee for the comments leading to improvements of this paper.

\section{Complex and real character varieties}

\label{sec:real-char}

Let $F_{r}$ be a rank $r$ free group ($r\in \mathbb{N}$) and $\mathbf{G}$ be a complex reductive algebraic group
defined over $\mathbb{R}$. We also assume $\mathbf{G}$ irreducible.
The \emph{$\mathbf{G}$-representation
variety of $F_{r}$} is defined as $
\mathfrak{R}_{r}(\mathbf{G}):=\hom(F_{r},\mathbf{G}).$
There is a homeomorphism between $\mathfrak{R}_{r}(\mathbf{G})$ and $\mathbf{G}^{r}$, given by the evaluation map which is defined over $\mathbb{R}$,
 if $\mathfrak{R}_{r}(\mathbf{G})$
is endowed with the compact-open topology (as defined on a space of
maps, with $F_{r}$ given the discrete topology, and
$\mathbf{G}$ the Euclidean topology from some affine embedding
$\mathbf{G}\subset \mathbb{C}^n$, $n\in \mathbb{N}$), and $\mathbf{G}^{r}$
with the product topology.
As $\mathbf{G}$ is a smooth affine variety defined over $\mathbb{R}$,
$\mathfrak{R}_{r}(\mathbf{G})$ is also a smooth affine variety and
it is defined over $\mathbb{R}$.
Consider now the action of $\mathbf{G}$ on $\mathfrak{R}_{r}(\mathbf{G})$
by conjugation. This
defines an action of $\mathbf{G}$ on the algebra $\mathbb{C}[\mathfrak{R}_{r}(\mathbf{G})]$
of regular functions on $\mathfrak{R}_{r}(\mathbf{G})$. Let $\mathbb{C}[\mathfrak{R}_{r}(\mathbf{G})]^{\mathbf{G}}$
denote the subalgebra of $\mathbf{G}$-invariant functions.
Since $\mathbf{G}$ is reductive the affine categorical
quotient may be defined as $
\mathfrak{X}_{r}(\mathbf{G}):=\mathfrak{R}_{r}(\mathbf{G})/\!/\mathbf{G}=\mathrm{Spec}_{\max}(\mathbb{C}[\mathfrak{R}_{r}(\mathbf{G})]^{\mathbf{G}})$.
This is a singular affine variety (irreducible and normal, since $\mathbf{G}^r$ is smooth and irreducible\footnote{
In Subsection 2.2 of \cite{Casimiro-Florentino-Lawton-Oliveira:2014} it was accidentally stated that $\X_r(\bG)$ is not necessarily irreducible; this is generally the case when $F_r$ is replaced by a finitely generated group $\Gamma$.}), whose points correspond to the Zariski
closures of the orbits. Since $\X_r(\bG)$ is an affine variety, it is a subset of an affine space, and inherits the Euclidean topology. With respect to this topology, in \cite{Florentino-Lawton:2013b}, it is shown that $\X_r(\bG)$ is homeomorphic to the conjugation orbit space of closed orbits (called the {\it polystable quotient}). $\mathfrak{X}_{r}(\mathbf{G})$, together with that topology, is called the \emph{$\mathbf{G}$-character variety}.

Let us define the conditions on a real Lie group $G$, for which
our results will apply.
\begin{defn}
\label{def:condforG} Let $K$ be a compact
Lie group. We say that $G$ is a {\em real $K$-reductive Lie group} if the following conditions hold:
(1) $K$ is a maximal compact subgroup of $G$;
(2) there exists a complex reductive algebraic group $\mathbf{G}$, defined over $\R$, such that
$\mathbf{G}(\R)_{0}\subseteq G\subseteq\mathbf{G}(\R),$
 where $\mathbf{G}(\R)$ denotes the real algebraic group of $\R$-points
of $\mathbf{G}$, and $\mathbf{G}(\R)_{0}$ its identity component (in the Euclidean topology);
(3) $G$ is Zariski dense in $\mathbf{G}$.
\end{defn}

We note that, if $G\neq\mathbf{G}(\R)$, then $G$ is not
necessarily an algebraic group (consider for example $G=\GL(n,\R)_{0}$).
When $K$ is understood, we often simply call $G$ a real reductive Lie group.
All classical real matrix groups, as well as all complex reductive Lie groups,
are in this setting. As an example which is
not under the conditions of Definition \ref{def:condforG}, we can
consider $\widetilde{\SL(n,\R)}$, the universal covering group of
$\SL(n,\R)$. For $n\geq 2$ it is not a matrix group, and so does not satisfy our definition, since all real reductive $K$-groups are linear.

%\vspace{0.5cm}

As above, let $K$ be a compact Lie group, and $G$ be a real $K$-reductive Lie group.
In like fashion, we define the \emph{$G$-representation variety of $F_{r}$}:
$
\mathfrak{R}_{r}(G):=\hom(F_{r},G).
$
Again, $\mathfrak{R}_{r}(G)$ is homeomorphic to $G^{r}$.
Similarly, as a set, we define
$\mathfrak{X}_{r}(G):=\mathfrak{R}_{r}(G)/\!/G
$
to be the set of closed orbits under the conjugation action of $G$
on $\mathfrak{R}_{r}(G)$. We give $\mathfrak{X}_{r}(G)$ the quotient topology on the subspace of points with closed orbits in $\mathfrak{R}_{r}(G)$.   This quotient coincides with the one considered by Richardson-Slodowy in \cite[Section 7]{rich-slod:1990}, and it is a non-trivial result in \cite{rich-slod:1990} that this quotient is always Hausdorff.  It is likewise called the \emph{$G$-character variety of $F_r$} even though it may not even be a semi-algebraic set.  However, it is an affine real semi-algebraic set when G is real algebraic.
For $K$ a compact Lie group, with its usual topology, we also define the space $\mathfrak{X}_{r}(K):=\hom(F_{r},K)/K\cong K^{r}/K,$ called the \emph{$K$-character variety} of $F_{r}$.  Since the $K$-orbits are always closed and $K$ is real algebraic, this construction is a special case of the construction above.  So it is Hausdorff and  can be identified with a semi-algebraic subset of $\mathbb{R}^{d}$, for some $d$. Moreover, it is compact, being the compact quotient of a compact space.

\section{Cartan decomposition and deformation to the maximal compact}

\label{sec:polar-dec}
Let $\lieg$ denote the Lie algebra of $G$, and $\liegc$
the Lie algebra of $\mathbf{G}$. We will fix a Cartan involution of $\liegc$
which restricts to a Cartan involution, $\theta$, of $\lieg$. This choice allows for a Cartan decomposition of those Lie algebras. In paricular, $\lieg=\liek\oplus\liep$ with $\theta|_{\liek}=1$ and $\theta|_{\liep}=-1$. Furthemore, $\liek$ is the Lie algebra of a maximal compact subgroup, $K$, of $G$. The Cartan involution of $\lieg$ lifts
to a Lie group involution $
\Theta:G\to G$
 whose differential is $\theta$ and such that $K=\Fix(\Theta)=\{g\in G:\, \Theta(g)=g\}$.
 The multiplication map
 provides a diffeomorphism $G\simeq K\times\exp(\liep)$. In particular, the exponential is injective on $\liep$.
 If we write $g=k\exp(X)$, for some $k\in K$ and $X\in\liep$, then
$
\Theta(g)^{-1}g=\exp(2X).$
So define $
(\Theta(g)^{-1}g)^{t}:=\exp\left(2tX\right),$
for any real parameter $t$.

\begin{prop}
\label{sdr}The map $H:[0,1]\times G\to G$, $H(t,g)=f_{t}(g):=g(\Theta(g)^{-1}g)^{-t/2}$ is
a strong deformation retraction from $G$ to $K$, and for each $t$,
$f_{t}$ is $K$-equivariant with respect to the action of
conjugation of $K$ in $G$. \end{prop}

By Proposition \ref{sdr}, there is a $K$-equivariant strong deformation
retraction from $G$ to $K$, so there is a $K$-equivariant strong
deformation retraction from $G^{r}$ onto $K^{r}$ with respect to
the diagonal action of $K$. This immediately implies:
\begin{cor}
\label{sdr1} Let $K$ be a compact Lie group and $G$ be a real $K$-reductive
Lie group. Then $\mathfrak{X}_{r}(K)$ is a strong deformation retract
of $\mathfrak{R}_{r}(G)/K$.
\end{cor}

\section{Kempf-Ness set and deformation retraction for character varieties}

\label{sec:def-retract}As before, fix a compact Lie group $K$, and
a real $K$-reductive Lie group $G$. Suppose that $G$ acts linearly
on a complex vector space $\mathbb{V}$, equipped with a Hermitian
inner product $\la\,,\ran$. Without loss of generality we can
assume that $\la\,,\ran$ is $K$-invariant, by averaging.
\begin{defn}
A vector $X\in\VV$ is a \emph{minimal vector} for the action of $G$
on $\VV$ if \newline
 $\|X\|\leq\|g\cdot X\|,$
for every $g\in G$, where $||\cdot||$ is the norm corresponding
to $\la\,,\ran$.
\end{defn}

Let $\mathcal{KN}_{G}=\mathcal{KN}(G,\VV)$
denote the set of minimal vectors. $\mathcal{KN}_{G}$ is known as
the \emph{Kempf-Ness set} in $\VV$ with respect to the action of
$G$. It is a closed algebraic set in $\mathbb{V}$.
The Kempf-Ness theory also works for closed $G$-subspaces. Indeed, let
$Y$ be an arbitrary closed $G$-invariant subspace of $\VV$, and %note: subvariety changed to subspace after resubmission
define $
\mathcal{KN}_{G}^{Y}:=\mathcal{KN}_{G}\cap Y.$
The next theorem is proved in \cite[Proposition 7.4, Theorems 7.6, 7.7 and 9.1]{rich-slod:1990}.
\begin{thm}
\label{thm:RS} The quotient $Y\quot G$ is a closed Hausdorff space and if, in addition, $Y$ is a real algebraic subset of $\VV$, it is also homeomorphic to a closed semi-algebraic set in some $\mathbb{R}^{d}$. Moreover, there is a $K$-equivariant deformation retraction of $Y$ onto $\mathcal{KN}_{G}^{Y}$.
\end{thm}

To apply the Kempf-Ness theorem to our situation, we need to embed
the $G$-invariant closed set $Y=\mathfrak{R}_{r}(G)=\hom(F_{r},G)\cong G^r$
in a complex vector space $\mathbb{V}$, as follows.
As $\mathbf{G}$ is a complex reductive algebraic group, $\mathbf{G}$ and $G$ can embedded in some $ \mathrm{GL}(n,\mathbb{C})$ and $ \mathrm{GL}(n,\mathbb{R})$, respectively. In this case, the Cartan involution is given by $\theta(A)=-A^{t}$, so that $\Theta(g)=(g^{-1})^{t}$. From now on, we will assume this situation.
We can obtain the embedding of $K^{r}$ ($r\in\mathbb{N}$)
into the vector space given by the product of the spaces of all $n$-square
complex matrices, which we denote by $\VV$.  Precisely, $\VV:=\mathfrak{gl}(n,\C)^{r}\cong\C^{rn^{2}}$.
The adjoint representation of $\GL(n,\C)$ in $\mathfrak{gl}(n,\C)$
restricts to a representation $
G\to\Aut(\VV)$
given by $
g\cdot(X_{1},\ldots,X_{r})=(gX_{1}g^{-1},\ldots,gX_{r}g^{-1}),\hspace{0.2cm}g\in G,\ X_{i}\in\mathfrak{gl}(n,\C).$
Moreover, this yields a representation $
\lieg\to\End(\VV)$
of the Lie algebra $\lieg$ of $G$ in $\VV$ given by the Lie brackets:
$
A\cdot(X_{1},\ldots,X_{r})=(AX_{1}-X_{1}A,\ldots,AX_{r}-X_{r}A)=([A,X_{1}],\ldots,[A,X_{r}])$
for every $A\in\lieg$ and $X_{i}\in\mathfrak{gl}(n,\C)$. In what
follows, the context will be clear enough to distinguish the notations
of the above representations.
We choose an inner product $\la\,,\ran$ in $\mathfrak{gl}(n,\C)$
which is $K$-invariant, under the restriction of the representation
$\GL(n,\C)\to\Aut(\mathfrak{gl}(n,\C))$ to $K$. From this we obtain
an inner product on $\VV$, $K$-invariant by the corresponding diagonal
action of $K$: $
\la(X_{1},\ldots,X_{r}),(Y_{1},\ldots,Y_{r})\ran=\sum_{i=1}^{r}\la X_{i},Y_{i}\ran$
for $X_{i},Y_{j}\in\mathfrak{gl}(n,\C)$. In $\mathfrak{gl}(n,\C)$,
$\la\,,\,\ran$ can be given explicitly by $\la A,B\ran=\tr(A^{*}B)$.
So by Theorem \ref{thm:RS}, Proposition \ref{sdr} and Corollary \ref{sdr1} we have the following theorem:
\begin{thm}
\label{main} The spaces $\mathfrak{X}_{r}(G)$ and $\mathfrak{X}_{r}(K)$
have the same homotopy type. In particular, the homotopy type of $\mathfrak{X}_{r}(G)$ depends only
on the maximal compact subgroup $K$ of $G$.\end{thm}
In our setting, the Kempf-Ness can be explicitly described as the closed set given by:
\begin{prop}
\label{pro:general-KN} $
\mathcal{KN}_{G}^{Y}=\left\{(g_{1},\cdots,g_{r})\in G^{r}:\ \sum_{i=1}^{r}g_{i}^{*}g_{i}=\sum_{i=1}^{r}g_{i}g_{i}^{*}\right\}.$
In particular, we have the inclusion $K^{r}\cong\hom(F_{r},K)\subset\mathcal{KN}_{G}^{Y}$.
\end{prop}
For $G$ algebraic, there is a natural inclusion of finite CW-complexes $\mathfrak{X}_{r}(K)\subset\mathfrak{X}_{r}(G)$ (see Lemma 4.9 of \cite{Casimiro-Florentino-Lawton-Oliveira:2014}). Using this result,
 Proposition \ref{pro:general-KN},  Theorem \ref{thm:RS} and Whitehead's Theorem we achieve the main result: (see Section 6 of \cite{Casimiro-Florentino-Lawton-Oliveira:2014} for an example)
\begin{thm}\label{thm:maintheorem}
There is a strong deformation retraction from $\mathfrak{X}_{r}(G)$
to $\mathfrak{X}_{r}(K)$.
\end{thm}

% Lista de bibliografia

\end{document}